\documentclass [12 pt, leqno] {amsart}
 
\usepackage {amssymb}
\usepackage {amsthm}
\usepackage {amscd}

\usepackage[all]{xy}
\usepackage {hyperref} 

\begin {document}

\theoremstyle{plain}

\newtheorem{lemma}[subsection]{Lemma}

\newtheorem{thm}[subsection]{Theorem}
\newtheorem*{thm*}{Theorem}

\theoremstyle{definition}

\theoremstyle{remark}
\newtheorem{remark}[subsection]{Remark}

\numberwithin{equation}{subsection}

\newcommand{\eq}[2]{\begin{equation}\label{#1}  #2 \end{equation}}

\newcommand{\ml}[2]{\begin{multline}\label{#1}  #2 \end{multline}}

\newcommand{\nnal}[1]{\begin{align*} #1 \end{align*}}
\newcommand{\nneq}[1]{\begin{equation} \nonumber #1 \end{equation}}

\newcommand{\nnml}[1]{\begin{multline}\nonumber #1 \end{multline}}


\newcommand{\arir}{\ar@{^{(}->}}
\newcommand{\aril}{\ar@{_{(}->}}

\newcommand{\are}{\ar@{>>}}

\newcommand{\xr}[1] {\xrightarrow{#1}}



\newcommand{\CH}{{\rm CH}}
\newcommand{\Gr}{{\rm Gr}}

\renewcommand{\H}{{\rm H}}
\newcommand{\id}{{\rm id}}


\newcommand{\Z} {\mathbb{Z}}
\newcommand{\Q} {\mathbb{Q}}

\renewcommand{\Gr}{\mathbb{G}{\rm r}}


\renewcommand{\P}{\mathbb{P}}

\renewcommand{\O}{\mathcal{O}}


\newcommand{\Sym}[1]{{\rm Sym}^{#1}}
\newcommand{\Sing}{{\rm Sing}}
\renewcommand{\H}{\mathcal{H}}
\newcommand{\rk}{{\rm rk}}
\newcommand{\Cor}{{\rm Cor}}

\title[Motives of hypersurfaces of small degree]{Motives of hypersurfaces of very small degree}
\author{Andre Chatzistamatiou}
\email{Andre.Chatzistamatiou@ens.fr}

\begin{abstract}
We study the Chow motive (with rational coefficients) of a hypersurface $X$ in the projective space
by using the variety $F(X)$ of $l$-dimensional planes contained in $X$.
If the degree of $X$ is sufficiently small we show that the
primitive part of the motive of $X$ is the tensor product of a direct summand in the motive of a suitable complete 
intersection in $F(X)$ and the $l$-th twist $\Q(-l)$ of the Lefschetz motive.   
\end{abstract}

\maketitle

\section*{Introduction} 

Let $X$ be a smooth hypersurface of degree $d$ in the projective space $\P^n_k$ over a field $k$. 
In this paper we study the Chow motive (with rational coefficients) of $X$ provided that $d$ is sufficiently small. 

Roitman has shown that the Chow group of zero-dimensional cycles of degree $0$ is a torsion group if $d\leq n$ \cite{R}.  
For higher dimensional cycles it is known \cite[Theorem~4.6]{ELV} that the Chow groups satisfy 
\nneq
{
\CH_{l'}(X)\otimes \Q = \CH_{l'}(\P^n)\otimes \Q = \Q \quad \text{for $0\leq l' \leq l-1$}
} 
if $n\geq {d+l-1 \choose l}$ and $d\geq 3$. 
The identity also holds if $X$ is covered by  $l$-dimensional planes \cite[Theorem~9.28]{V}, or more generally if
$X$ is a hyperplane section of a hypersurface $Y$ which is covered by  $l$-dimensional planes \cite{O}. 

Results on triviality of Chow groups give rise to a decomposition of the motive $\bar{X}$ associated with $X$. In our case, 
we get 
\nneq
{
\bar{X} \cong  M_D \otimes \Q(-l) \oplus \bigoplus_{i=0}^{l-1} \Q(-i), 
}
where $\Q(-1)$ is the Lefschetz motive and $M_D$ is a direct summand of the motive $\bar{D}$ of some variety $D$. 
Our purpose is to describe $M_D$.

In order to state the theorem we need the following notation. An $l$-dimensional plane $E$ in $\P^n$ 
is called osculating plane if the intersection $E\cap X$ in an $(l-1)$-dimensional plane or if $E$ is contained in $X$. 
We say that $X$ has sufficiently many osculating planes if there exists an osculating plane through every closed 
point of $X$ (the planes may be defined over a field extension of $k$).    

\begin{thm*}[see Theorem \ref{mainthm}] 
Let $n,d,l$ be numbers such that a general hypersurface of degree $d$ in the projective space $\P^n$ has sufficiently 
many osculating $l$-dimensional planes.
Let $X\subset \P^n$ be a smooth hypersurface of degree $d$ such that the Fano variety $F_{l}(X)$ of $l$-dimensional planes 
contained in $X$ is smooth and has the expected dimension. Furthermore, let $HF_{l}(X)\subset F_{l}(X)$ be a smooth complete
intersection of hyperplanes (in the Pl\"ucker embedding) with $\dim HF_{l}(X)=n-2l-1$. 
Then there is an isomorphism in the category of Chow motives with rational coefficients:
$$
\bar{X} \cong M_{HF_{l}(X)} \otimes \Q(-l) \oplus \bigoplus_{i=0}^{n-1} \Q(-i)
$$ 
and $M_{HF_{l}(X)}$ is a direct summand in the motive of $HF_{l}(X)$.
\end{thm*} 

The conditions on $n,d,l$ hold if $n\geq {l+d-1 \choose l}+l-1$. 

Let us sketch the idea of the proof. We consider the family of planes over $HF_{l}(X)$:
$$
\Xi = \{(x,E) \in X \times HF_{l}(X); x\in E\} \subset X \times HF_{l}(X).
$$
The cycle $\Xi$ defines a correspondence $\phi_1: HF_{l}(X)\otimes \Q(-l) \xr{} X$ resp. 
$\phi_2: X \xr{} HF_{l}(X)\otimes \Q(-l)$. The composite $\phi_1 \circ \phi_2$ is the cycle    
$$
Z_X = \{ (x,y)\in X\times X; x,y\in E, E\in  HF_{l}(X) \} 
$$
in $\CH^{n-1}(X\times X)$. The most important step is to show that 
\eq{introdec}
{
Z_X = \imath^*(a) + m \cdot \Delta_X,
}
for some $a\in \CH^{n-1}(\P^n\times \P^n)$, some nonzero integer $m$, the inclusion $\imath:X^2 \hookrightarrow (\P^n)^2$ and  
the diagonal $\Delta_X$. 
In order to prove \ref{introdec} we introduce the doubled incidence variety of degree $d$ hypersurfaces together 
with two points:  
$$
\Sigma = \{(x,y,Y) \in (\P^n)^2 \times \P(\Sym{d}(k^{n+1})); x,y \in Y \}.
$$
The projection $p:\Sigma \xr{} (\P^n)^2$ is a projective bundle over $(\P^n)^2-\Delta_{\P^n}$ and $\Delta_{\P^n}$, 
so that $\CH^{n-1}(\Sigma)$ can be calculated by using the projective bundle formula and the localization sequence.
The cycle $p^{-1}(\Delta_{\P^n})\in \CH^{n-1}(\Sigma)$ maps to the diagonal $\Delta_X$ by the pullback map $\jmath^*$ of the 
inclusion $\jmath: X\times X\xr{} \Sigma$. One defines a cycle $Z\in \CH^{n-1}(\Sigma)$ with $\jmath^*(Z)=Z_X$, and by comparing
$Z$ with $p^{-1}(\Delta_{\P^n})$ we obtain \ref{introdec} after applying $\jmath^*$. 

\subsection*{Acknowledgments}
I thank H. Esnault for introducing me to this subject and for her interest in my work. I am grateful to Y. Andr\'e for valuable discussions. 
This paper is written during a stay at the Ecole normale sup\'erieure which is supported by a fellowship within the Post-Doc program of the Deutsche Forschungsgemeinschaft (DFG).
I thank the Ecole normale sup\'erieure for its hospitality.

\section{Cycles in the doubled incidence variety}

\subsection{}
Let $k$ be a field and let $\P(\Sym{d}(k^{n+1}))$ be the hypersurfaces of degree $d$ in the projective space $\P^n_k$. 
We denote by $\Sigma$ the doubled incidence variety 
$$
\Sigma = \{(x,y,X) \in (\P^n)^2 \times \P(\Sym{d}(k^{n+1})); x,y \in X \}.
$$
Let $p:\Sigma \xr{} (\P^n)^2$ be the projection, we define $\Sigma_1:=p^{-1}(\Delta_{\P^n})$ and $\Sigma_0:=\Sigma - \Sigma_1$.
In the diagram 
$$
\xymatrix
{
\Sigma_0 \ar[r] \ar[d]
&
\Sigma \ar[d]_{p}
&
\Sigma_1 \ar[l] \ar[d]
\\
(\P^n)^2 - \Delta_{\P^n} \ar[r]
&
(\P^n)^2
&
\Delta_{\P^n} \ar[l]
}
$$
the varieties $\Sigma_0$ resp. $\Sigma_1$ are projective bundle with fiber dimension $N-2$ resp. $N-1$, 
where $N=\dim \P(\Sym{d}(k^{n+1}))$. It is easy to see that $\Sing(\Sigma)=\{(x,x,X); x \in \Sing(X)\}$. 
The singular locus is a projective bundle over $\Delta_{\P^n}$ with fiber dimension $N-1-n$.

\subsection{} 
There is an exact sequence 
\eq{locseq}
{
\CH^0(\Sigma_1) \xr{} \CH^{n-1}(\Sigma) \xr{} \CH^{n-1}(\Sigma_0) \xr{} 0
} 
from the localization sequence of Chow groups ($\CH^i$ denotes the group of $i$-codimensional cycles
modulo rational equivalence). Moreover, there is a natural splitting defined as follows.  

Let \small$\phi:\Sigma\xr{} \P(\Sym{d}(k^{n+1}))$ \normalsize be the projection and set \small$c=\phi^*(c_1(\O(1)))$. \normalsize
For the other projection \small$p: \Sigma \xr{} (\P^n)^2$  \normalsize we may define $p^*$ to be the composite $\epsilon^* \circ {\rm pr}^*$ where
\small$\epsilon: \Sigma \xr{} (\P^n)^2 \times \P(\Sym{d}(k^{n+1}))$  \normalsize is the regular embedding and ${\rm pr}$ is the projection to $(\P^n)^2$. 
From the projective bundle formula we see that 
$$
\bigoplus_{i=0}^{n-1} c^{n-1-i}\cdot p^*\CH^i((\P^n)^2) \subset \CH^{n-1}(\Sigma)
$$
splits the sequence \ref{locseq}, so that every class $Z$ in $\CH^{n-1}(\Sigma)$ can be written as
$$
Z= \sum_{i=0}^{n-1} c^{n-1-i}\cdot p^*(a_i) + m \cdot [\Sigma_1]
$$
with $a_i\in \CH^{i}((\P^n)^2)$ and $m\in \Z$.

\subsection{}
Let $\H$ be a smooth, connected projective  $k$-scheme and $\Xi\subset \H \times \P^n$ be a family of $\kappa$ dimensional subschemes
of $\P^n$, flat over $\H$. We assume that $\kappa\geq 1$ and:
\begin{itemize}
\item[(A)] The sheaf ${\rm pr}_{1*}\left( \O_{\Xi} \otimes {\rm pr}^*_2\O_{\P^n}(d) \right)$ on $\H$ is locally free and the 
natural map 
\eq{SymOd}
{
\Sym{d}((k^{n+1})^{\vee})\otimes \O \xr{} {\rm pr}_{1*}\left( \O_{\Xi} \otimes {\rm pr}^*_2\O_{\P^n}(d) \right)
}
is surjective.
\end{itemize}

We denote by $E$ the kernel of \ref{SymOd}. It is convenient to write $QE$ for ${\rm pr}_{1*}\left( \O_{\Xi} \otimes {\rm pr}^*_2\O_{\P^n}(d) \right)$.
In the following commutative diagram we fix the notation for the various maps
$$
\xymatrix
{
\Xi \times_{\H}\Xi \times_{\H} \P(E^{\vee}) \ar[r]^-{f_{\Sigma}} \ar[d]_-{f_{\P(E^{\vee})}}
&
\Sigma \ar[r]^-{p} \ar[d]^-{\phi}
&
(\P^n)^2 
\\
\P(E^{\vee}) \ar[r]_-{\psi} \ar[d]_-{f_{\H}}
&
\P(\Sym{d}(k^{n+1}))
\\
\H.
}
$$
Define 
\eq{defZ}
{
Z:= f_{\Sigma*}f_{\P(E^{\vee})}^*f_{\H}^*:\CH^{*-e}(\H) \xr{} \CH^*(\Sigma),
} 
we will be mainly interested in cycles $Z(a)\in \CH^{n-1}(\Sigma)$ and their pullback to $X\times X\subset \Sigma$ for 
a hypersurface $X$.

The cycle $\Xi\in \CH^{n-\kappa}(\H \times \P^n)$ has a unique representation 
\eq{Xi}
{
[\Xi] = \sum_{i=0}^{n-\kappa} \xi_{n-\kappa-i} \otimes H^i  
}  
with $\xi_{j}\in \CH^j(\H)$ and $H=c_1(\O_{\P^n}(1))$.

\begin{lemma} \label{lemma:ai}
For $a\in \CH^{n-1-e}(\H)$ let $Z(a)_{\mid \Sigma_0}=\sum_{i=0}^{n-1} c^{n-1-i} \cdot p^*(a_i)$ be the pullback 
of $Z(a)$ to $\Sigma_0$. The classes $a_i$ can be computed as follows:
$$
\sum_{i} a_i = \frac{\sum_{0\leq s,t \leq n-\kappa} b_{s,t}H^s\otimes H^t}{(1+d\otimes H)(1+H\otimes d)}
$$ 
in $\CH^*((\P^n)^2-\Delta_{\P^n})$, and where 
$$
b_{s,t}=\int_{\H} \xi_{n-\kappa-s} \cdot \xi_{n-\kappa-t}\cdot c_{\rk(QE)-n+s+t-1}(QE)\cdot a
$$
\begin{proof}
On $(\P^n)^2-\Delta_{\P^n}$ the evaluation morphism 
\eq{ev}
{
\Sym{d}((k^{n+1})^{\vee})\otimes \O \xr{} {\rm pr}_1^*\O_{\P^n}(d)\oplus {\rm pr}_2^*\O_{\P^n}(d)
}
is surjective, let $G$ be the kernel. We have $\Sigma_0=\P(G^{\vee})$ and $c_{\mid \Sigma_0}=c_1(\O_{\P(G^{\vee})}(1))$, so that
\eq{SegreG}
{
p_*(\frac{1}{1-c} Z(a)_{\mid \Sigma_0})= \sum_{j\geq 0} s_j(G) \cdot \sum_{i=0}^{n-1} a_i,
}
where $s_j(G)$ are the Segre classes. Since $G$ is the kernel of \ref{ev} we see 
\eq{SegreG2}
{
(1+d\cdot {\rm pr}_1^*c_1(\O_{\P^n}(1)))(1+d\cdot {\rm pr}_2^*c_1(\O_{\P^n}(1))) \sum_{j\geq 0} s_j(G) = 1.
}
Define 
\eq{Ti}
{
T_i(a) = p_*(c^{N-n-1+i} Z(a))
}
and let $\jmath:(\P^n)^2-\Delta_{\P^n}\subset (\P^n)^2$ be the open immersion. It follows from \ref{SegreG} and \ref{SegreG2} that
\eq{Tiai}
{
\sum_{i=0}^{n-1} a_i = \frac{ \sum_{i\geq 0} \jmath^*T_i(a)}{(1+d\otimes H)(1+H\otimes d)}.
} 
Let us now compute the $T_i$:
\nnal
{
p_*(c^{N-n-1+i} Z(a)) &= p_* f_{\Sigma*}((\phi\circ f_{\Sigma})^* c_1(\O(1))^{N-n-1+i} \cdot (f_{\H}\circ f_{\P(E^{\vee})})^*(a)) \\
                      &= p_* f_{\Sigma*}f_{\P(E^{\vee})}^*( \psi^* c_1(\O(1))^{N-n-1+i} \cdot  f_{\H}^*(a) ) 
}
The map \small$T=p_* f_{\Sigma*}f_{\P(E^{\vee})}^*$\normalsize is given by the correspondence \small$\P(E^{\vee}) \times_{H}\Xi \times_{\H} \Xi$ \normalsize in 
$\P(E^{\vee})\times (\P^n)^2$, and using \ref{Xi} we see
$$
[\Xi \times_{\H} \Xi] = \sum_{0\leq s,t \leq n-\kappa} (\xi_{n-\kappa-s}\cdot \xi_{n-\kappa-t}) \otimes H^s\otimes H^t
$$
in $\CH^*(\H\times (\P^n)^2)$. So that the coefficient of $H^s\otimes H^t$ in $T_i(a)$ is
\ml{coefficientTi}
{
\int_{\P(E^{\vee})}  f_{\H}^*\xi_{n-\kappa-s}\cdot f_{\H}^*\xi_{n-\kappa-t} \cdot\psi^* c_1(\O(1))^{N-n-1+i} \cdot  f_{\H}^*(a) =\\
  \int_{\H}  \xi_{n-\kappa-s}\cdot \xi_{n-\kappa-t} \cdot s_{N-n-\rk(E)+i}(E) \cdot a.  
}
If $Z(a)\in \CH^{n-1}(\Sigma)$ then $T_i(a)\in \CH^{i}((\P^n)^2)$ by definition, so that the coefficient of $H^s\otimes H^t$ in 
$T_i(a)$ vanishes if $s+t\neq i$. The identity $s_{N-n-\rk(E)+i}(E)=c_{\rk(QE)-n-1+i}(QE)$ completes the prove.
\end{proof}
\end{lemma}

The next Lemma computes the pullback of $Z(a)$ to $X\times X$ for a hypersurface $X$. We write $\imath_{\Sigma}$ for the
inclusion $X\times X \xr{} \Sigma$ and $\imath_{(\P^n)^2}$ for the inclusion $X\times X \xr{}(\P^n)^2$. Note that both inclusion
are locally complete intersection, thus the pullback is well-defined.   

\begin{lemma} \label{lemma:m}
For $Z(a)\in \CH^{n-1}(\Sigma)$ and $Z(a)_{\mid \Sigma_0}=\sum_{i=0}^{n-1} c^{n-1-i} \cdot p^*(a_i)$, we have
$$
\imath_{\Sigma}^* Z(a) = \imath_{(\P^n)^2}^*(a_{n-1}) - m \cdot \Delta_X,
$$
where 
$$
m=d \cdot \sum_{j=\kappa-1}^{n-1} (-d)^j  \int_{\H} \xi_{n-\kappa} \cdot \xi_{j-\kappa+1}\cdot c_{\rk(QE)-2-j}(QE)\cdot a.
$$
\begin{proof}
We know that 
\eq{repr}
{
Z(a)=\sum_{i=0}^{n-1} c^{n-1-i} \cdot p^*(a_i) - m\cdot [\Sigma_1]
}
for some $m$. The line bundle $\imath_{\Sigma}^*\phi^*\O(1)$ is trivial and $\imath_{\Sigma}^*[\Sigma_1]=\Delta_X$, therefore
\eq{onXX}
{
\imath_{\Sigma}^* Z(a)=\imath_{(\P^n)^2}^* a_{n-1} - m \cdot \Delta_X.
}
We claim that 
\eq{claim1}
{
\imath_{(\P^n)^2 *}\imath_{\Sigma}^* \beta = p_*(c^N \cdot \beta)
}
for every class $\beta \in \CH^*(\Sigma)$. This follows from the following fact. If $g:D\subset Y$ 
is a Cartier divisor on $Y$ and $L$ the associated line bundle, then $g_*g^*(\beta)=c_1(L)\cdot \beta$.

Claim \ref{claim1} implies $\imath_{(\P^n)^2 *}\imath_{\Sigma}^* Z(a)=T_{n+1}(a)$ (see \ref{Ti}). The coefficient 
of $H\otimes H^n$ in $T_{n+1}(a)$ is computed in \ref{coefficientTi} and vanishes for trivial reasons (if $\kappa\geq 1$).
Using Lemma \ref{lemma:ai}, we see that the coefficient $\gamma$ of $1\otimes H^{n-1}$ in $a_{n-1}$ is 
$$
\gamma = \sum_{j=\kappa-1}^{n-1} (-d)^j \cdot b_{0,n-1-j}.
$$
By applying $\imath_{(\P^n)^2*}$ to \ref{onXX} and using \small$\imath_{(\P^n)^2*} \imath_{(\P^n)^2}^* a_{n-1} = d^2 (H\otimes H)\cdot a_{n-1}$ \normalsize and 
\small$\imath_{(\P^n)^2*} \Delta_X= d\cdot \sum_{i\geq 0}^{n-1}H^{i+1} \otimes H^{n-i}$ \normalsize it follows that 
$
m= d \cdot \gamma,
$
as claimed.
\end{proof}
\end{lemma}

\section{Motives of hypersurfaces and their Fano varieties}

\subsection{} In the following we work with the Grassmannian of $\kappa$-planes $\H=\Gr_{\kappa}$  in projective space 
and $\Xi\subset \Gr_{\kappa}\times \P^n$ the universal family. We denote by $V$ resp. $QV$ the tautological 
bundle $V\subset \O_{\Gr}^{n+1}$ resp. the quotient $\O_{\Gr}^{n+1}/V$. The family $\Xi$ is the projective bundle $\Xi=\P(V^{\vee})$
and it is easy to see that 
$$
[\Xi] = \sum_{i=0}^{n-\kappa} c_{n-\kappa-i}(QV)\otimes H^i
$$     
in $\CH^{n-\kappa}(\Gr\times \P^n)$. Furthermore we have 
$QE={\rm pr}_{1*}( \O_{\Xi} \otimes {\rm pr}^*_2\O_{\P^n}(d))=\Sym{d}(V^{\vee})$.

\subsection{} We will be interested in cycles $Z(c_1(V^{\vee})^s)\in \CH^{n-1}(\Sigma)$, for $s\geq 0$, (notation as in \ref{defZ}).
By counting dimensions we see that 
\eq{s}
{
s=\dim \Xi - \rk(QE) - (n-1) = \kappa(n-\kappa)-{d+\kappa \choose \kappa} + \kappa + 1
} 
Let us consider the variety 
\ml{oscplanes}
{
\{ (x,E_{\kappa-1},E_{\kappa},X) \in \P^n \times \Gr_{\kappa-1} \times \Gr_{\kappa} \times  \P(\Sym{d}(k^{n+1})); \\
x\in E_{\kappa-1}\subset E_{\kappa},  \text{$E_{\kappa}\cap X = E_{\kappa-1}$ or $E_{\kappa}\subset X$}  \}.
}
More formally, this variety is defined as follows. On $\Xi=\P(V^{\vee})$ there is an exact sequence of vector bundles
\eq{exseqV1V}
{
0 \xr{} V_1^{\vee} \xr{} V^{\vee} \xr{} \O_{\P(V^{\vee})}(1) \xr{} 0,
} 
and the points of $\P(V_1)$ are $\{(x,E_{\kappa-1},E_{\kappa}); x\in E_{\kappa-1}\subset E_{\kappa}\}$. Since $\O_{\P(V_1)}(-1)\subset 
V_1^{\vee} \subset V^{\vee}$ we can define $G$ to be the kernel of 
$$
\Sym{d}((k^{n+1})^{\vee}) \xr{} \Sym{d}(V^{\vee})/\O_{\P(V_1)}(-d).
$$
Then $\P(G^{\vee})$ is the variety \ref{oscplanes}. 

The following condition will imply that the diagonal $\Delta_X$, for 
a hypersurface $X$, can be written in terms of the pullback of $Z(c_1(V^{\vee})^s)$ to $X\times X$ 
(i.e. $m\neq 0$ in Lemma \ref{lemma:m}).
\begin{itemize}
\item[(B)] The following map is surjective:
\nnal
{
 &\P(G^{\vee}) \xr{} \{(x,X)\in \P^n \times \P(\Sym{d}(k^{n+1})); x \in X \} \\
 &(x,E_{\kappa-1},E_{\kappa},X)  \mapsto (x,X).
}
\end{itemize} 
By counting dimensions we see that a necessary condition for (B) is $s\geq 0$
(with $s$ as in \ref{s}). If $d=2$ then $s\geq 0$ is not sufficient, the first example is $\kappa=3$ and $n=5$. 
In fact, (B) is equivalent to $n\geq 2\cdot \kappa$ if $d=2$, which can be checked by using the following Lemma.
However, we don't know any counter examples to 
\eq{equiv}
{
(B) \Leftrightarrow s\geq 0 
}  
for $d>2$. If $\kappa=1$ then \ref{equiv} is true. 
This also holds for $(n,d,\kappa)=(6,3,2),(8,4,2),(11,5,2),(9,3,3)$.
 
It is known \cite[Lemma~1.1+Lemma~4.2]{ELV} that (B) is true if $d\geq 3$ and  
$$
n-\kappa+1 \geq {\kappa-1+d \choose \kappa}. 
$$ 

\begin{lemma}\label{lemma:conditionBm}
Condition (B) holds if and only if 
$$
m=d \cdot \sum_{j=\kappa-1}^{n-1} (-d)^j  \int_{\H} \xi_{n-\kappa} \cdot \xi_{j-\kappa+1}\cdot c_{\rk(QE)-2-j}(QE)\cdot c_1(V^{\vee})^s
$$
is nonzero.
\begin{proof}
From the construction of $\P(G^{\vee})$ we have the maps
$$
\P(G^{\vee}) \xr{f} \P(V_1) \xr{g} \P(V^{\vee}) \xr{h} \Gr_{\kappa}.
$$
We claim that 
\ml{claim}
{
(h\circ g \circ f)_*( c_1(\O_{\P(V^{\vee})}(1))^{n-1} \cdot  c_1(\O_{\P(G^{\vee})}(1))^{N}) = \\ (-1)^{\kappa-1} d \cdot 
\sum_{j=\kappa-1}^{n-1} (-d)^j  \xi_{n-\kappa} \cdot \xi_{j-\kappa+1}\cdot c_{\rk(QE)-2-j}(QE)
}
($N=\dim \P(\Sym{d}(k^{n+1}))$). Indeed,  
\small 
\nnml
{
g_*f_*(c_1(\O_{\P(G^{\vee})}(1))^{N}) = g_*s_{\rk(QE)-1}(G) = g_*c_{\rk(QE)-1}(QE/\O_{\P(V_1)}(-d))=\\ 
g_*\sum_{i} d^i  c_{\rk(QE)-i-1}(QE) \cdot c_1(\O_{\P(V_1)}(1))^i = \sum_{i} d^i  c_{\rk(QE)-i-1}(QE) \cdot s_{i-\kappa+1}(V_1^{\vee}),
} \normalsize
and from \ref{exseqV1V} we obtain \small$s_j(V_1^{\vee})=s_j(V^{\vee})+c_1(\O_{\P(V^{\vee})}(1))\cdot s_{j-1}(V^{\vee})$ \normalsize 
for all $j$. Thus, 
\ml{claimfast}
{
(h\circ g \circ f)_*( c_1(\O_{\P(V^{\vee})}(1))^{n-1} \cdot  c_1(\O_{\P(G^{\vee})}(1))^{N}) = \\
     \sum_{i} d^i  c_{\rk(QE)-i-1}(QE) \cdot (s_{n-\kappa-1}(V)s_{i-\kappa+1}(V^{\vee})+s_{n-\kappa}(V)s_{i-\kappa}(V^{\vee})).
}

On $\P(V^{\vee})$ the natural morphism $QE\otimes \O_{\P(V^{\vee})}(-d)\xr{} \O$ is surjective, so that the top Chern class of $QE\otimes \O_{\P(V^{\vee})}(-d)$ vanishes:
\small
\nnml
{
h_*c_{\rk(QE)}(QE\otimes \O_{\P(V^{\vee})}(-d))=h_*\sum_i (-d)^i c_{\rk(QE)-i}(QE)\cdot c_1(\O_{\P(V^{\vee})}(1))^i= \\
\sum_i (-d)^i c_{\rk(QE)-i}(QE)s_{i-\kappa}(V) =(-1)^{\kappa}d\cdot \sum_i d^i c_{\rk(QE)-i-1}(QE)s_{i+1-\kappa}(V^{\vee})
}
\normalsize
is zero and together with \ref{claimfast} we obtain 
\nnml
{
(h\circ g \circ f)_*( c_1(\O_{\P(V^{\vee})}(1))^{n-1} \cdot  c_1(\O_{\P(G^{\vee})}(1))^{N}) = \\
  d \cdot \sum_{i} d^i  c_{\rk(QE)-i-2}(QE) \cdot s_{n-\kappa}(V)s_{i+1-\kappa}(V^{\vee}).
}
Then, $s_j(V^{\vee})=(-1)^js_j(V)=c_j(QV)=\xi_j$ proves the claim. 

Let $\pi$ be the map in condition (B). For a general closed point $(x,X)$ the irreducible components 
of $\pi^{-1}(x,X)$ map generically one to one to $\Gr_{\kappa}$. 
The class $c_1(V^{\vee})=c_1(\Lambda^{\kappa+1}V^{\vee})$ is the class of an ample line bundle and $s$ is the dimension of the generic
fiber of $\pi$, if $\pi$ is surjective. Thus, $\pi$ is surjective if and only if 
$\int_{\P(G^{\vee})}  c_1(V^{\vee})^s\cdot [\pi^{-1}(x,X)] \neq 0$. The class $c_1(\O_{\P(V^{\vee})}(1))^{n-1} \cdot  c_1(\O_{\P(G^{\vee})}(1))^{N}$
is equal to $\sum_i \pi^*(x_i,X_i)$ for some general points $(x_i,X_i)$ and \ref{claim} implies the Lemma. 
\end{proof}
\end{lemma}

\begin{remark}
From the proof of Lemma \ref{lemma:conditionBm} we obtain a geometric interpretation for $m$. Let $H\Gr_{\kappa}$ 
be a complete intersection of $s$ hyperplanes. For a general hypersurface $X$ and a general point $x\in X$ the number 
$m$ is up to the factor $(-1)^{\kappa-1}d$ the number of $\kappa$-planes in $H\Gr_{\kappa}$ through $x$ which meet $X$ in a 
$\kappa - 1$ plane with multiplicity $d$.    
\end{remark}

In the following Theorem we work with the category of (pure) Chow motives with rational coefficients (see \cite[Chapter~16]{F}) and $\Q(-1)$ denotes 
the Lefschetz motive.  

\begin{thm} \label{mainthm}
Let $n,d,\kappa$ be numbers satisfying (B). Let $X\subset \P^n$ be a smooth hypersurface of degree $d$ such that 
the Fano variety $F_{\kappa}(X)$ of $\kappa$-dimensional planes contained in $X$ has the expected dimension 
(which is $\dim F_{\kappa}(X)=(\kappa+1)(n-\kappa)-{d+\kappa \choose \kappa}$), and let $HF_{\kappa}(X)\subset F_{\kappa}(X)$ be a complete
intersection of hyperplanes (in the Pl\"ucker embedding) with $\dim HF_{\kappa}(X)=n-2\kappa-1$. Furthermore, let 
$\psi: \widetilde{HF_{\kappa}(X)} \xr{} HF_{\kappa}(X)$ be a generically finite surjective morphism which resolves the 
singularities of $HF_{\kappa}(X)$. Then there is an isomorphism in the category of Chow motives with rational coefficients:
$$
(X,\id_X) \cong (\widetilde{HF_{\kappa}(X)},P) \otimes \Q(-\kappa) \oplus \bigoplus_{i=0}^{n-1} \Q(-i)
$$ 
for a suitable projector $P$.
\begin{proof}
From Lemma \ref{lemma:m} and Lemma  \ref{lemma:conditionBm} we obtain 
\eq{decompdiag}
{
\Delta_X = -\frac{1}{m} \imath_{\Sigma}^* Z(c_1(V^{\vee})^s) + \frac{1}{m} \imath_{(\P^n)^2}^*(a_{n-1}).
}
In the cartesian diagram 
$$
\xymatrix
{
\Xi \times_{F_{\kappa}(X)}\Xi \ar[r]^-{f_{X\times X}} \ar[d]
&
X \times X \ar[d]^{\imath_{\Sigma}}
\\
\Xi \times_{\Gr_{\kappa}}\Xi \times_{\Gr_{\kappa}} \P(E^{\vee}) \ar[r]^-{f_{\Sigma}} 
&
\Sigma 
}
$$ 
pullback and pushforward commute, because $F_{\kappa}(X)$ has the expected dimension. Therefore 
\eq{firststep}
{
\imath_{\Sigma}^* Z(c_1(V^{\vee})^s) = f_{X\times X*}(c_1(V^{\vee})^s) = (f_{X\times X}\circ \eta)_*([\Xi \times_{HF_{\kappa}(X)}\Xi]),
}
where $\eta:HF_{\kappa}(X) \xr{} F_{\kappa}(X)$ is a complete intersection of $s$ hyperplanes 
(thus the dimension is $\dim HF_{\kappa}(X)=n-2\kappa-1$). It is more convenient to write $H=HF_{\kappa}(X)$ and 
$\tilde{H}=\widetilde{HF_{\kappa}(X)}$. 

There is a cycle $Y\in \CH_{n-2\kappa-1}(\tilde{H})\otimes \Q$ (i.e. $Y$ is a rational linear combination of 
connected components of $\tilde{H}$) such that $\psi_*(Y)=[H]$. 

Let $\phi_1\in \Cor(\tilde{H}\otimes  \Q(-\kappa),X)$ 
(resp. $\phi_2\in \Cor(X,\tilde{H}\otimes  \Q(-\kappa))$) be the correspondence defined by the cycle 
$\tilde{H}\times_{H} \Xi$ in $\tilde{H}\times X$ 
(resp. $[\Xi\times_{H}\tilde{H}]\cdot pr_{\tilde{H}}^*(Y)$ in $X \times \tilde{H}$, where
$pr_{\tilde{H}}:X\times  \tilde{H} \xr{}  \tilde{H}$ is the projection). 

We consider the commutative diagram 
$$
\xymatrix
{
&
\Xi\times_{\tilde{H}} \Xi \ar[dl] \ar[r]^{\psi'} \ar[dd]^{\cap}
&
\Xi\times_{H} \Xi \ar[dd]^{\cap} \ar[dr]^{f_{X\times X}\circ \eta}
&
\\
\tilde{H}
&
&
&
X\times X
\\
&
X\times \tilde{H}\times X \ar[ul]^{pr_{\tilde{H}}} \ar[r]^{\id \times \psi \times \id}
&
X\times H \times X \ar[ur]
&
}
$$
It is easy to see that $\phi_1\circ \phi_2=(f_{X\times X}\circ \eta\circ \psi')_*([\Xi\times_{\tilde{H}} \Xi]\cdot pr_{\tilde{H}}^*(Y))$, 
and $\psi'_*([\Xi\times_{\tilde{H}} \Xi]\cdot pr_{\tilde{H}}^*(Y))=[\Xi\times_{H} \Xi]$ together with \ref{firststep} yields 
\eq{secondstep}
{
\phi_1 \circ \phi_2 = Z(c_1(V^{\vee})^s).
}
 
If $H$ is the class of a hyperplane in $\P^n$ then we write $P_i$ for the 
pullback of $ \frac{1}{d} H^{n-1-i}\otimes H^{i}\in \CH^{n-1}(\P^n\times \P^n)\otimes \Q$ to $X\times X$. The correspondences 
$P_0,\dots,P_{n-1}$ are idempotent and orthogonal. We may write $a_{n-1}=\frac{m}{d}\sum_{i=0}^{n-1} \beta_i  H^{n-1-i}\otimes H^{i}$
and it follows from \ref{decompdiag} and \ref{secondstep} that
$$
\Delta_{X} - \sum_{i} P_i = -\frac{1}{m} \phi_1\circ \phi_2 + \sum_{i=0}^{n-1} (\beta_i-1)  P_i.
$$
Composition with $P_i$ shows that $\frac{1}{m} \phi_1\circ \phi_2\circ P_i =  (\beta_i-1)  P_i$ and 
\eq{thirdstep}
{
\Delta_{X} - \sum_{i} P_i = \phi_1 \circ (-\frac{1}{m} \phi_2 + \frac{1}{m} \sum_{i=0}^{n-1}  \phi_2\circ P_i).
}
Since \small $(X,P_i)\cong \Q(-i)$ \normalsize and \small$(X,\Delta_{X} - \sum_{i} P_i)\cong (\tilde{H}, -\frac{1}{m} \phi_2 \circ (\Delta_{X} - \sum_{i} P_i)\circ \phi_1)$ \normalsize
by \ref{thirdstep}, this proves the theorem. 
\end{proof}
\end{thm}

\end{document}